\documentclass{gtmon_a}
\pdfoutput=1

\proceedingstitle{Proceedings of the School and Conference in Algebraic
Topology (The Vietnam National University, Hanoi, 9-20 August 2004)}
\conferencestart{9 August 2004}
\conferenceend{20 August 2004}
\conferencename{School and Conference in Algebraic Topology}
\conferencelocation{Vietnam National University, Hanoi, Vietnam}

\editor{John Hubbuck}
\givenname{John}
\surname{Hubbuck}

\editor{Nguy\~\ecircumflex{}n H V H\uhorn{}ng}
\givenname{H\uhorn{}ng}
\surname{Nguy\~\ecircumflex{}n}

\editor{Lionel Schwartz}
\givenname{Lionel}
\surname{Schwartz}

\title{An algebraic introduction to the Steenrod algebra}

\author{Larry Smith}
\givenname{Larry}
\surname{Smith}
\address{AG-Invariantentheorie\\
Mittelweg 3\\
D 37133 Friedland\\
Federal Republic of Germany}
\email{larry@uni-math.gwdg.de}
\urladdr{}

\volumenumber{11}
\issuenumber{}
\publicationyear{2007}
\papernumber{15}
\startpage{327}
\endpage{348}

\doi{}
\MR{}
\Zbl{}

\arxivreference{}



\received{9 March 2005}
\revised{}
\accepted{20 June 2005}
\published{14 November 2007}
\publishedonline{14 November 2007}
\proposed{}
\seconded{}
\corresponding{}
\editor{}
\version{}

\makeatletter
\def\cnewtheorem#1[#2]#3{\newtheorem{#1}{#3}[section]
\expandafter\let\csname c@#1\endcsname\c@theorem}

\newtheorem{theorem}{Theorem}[section]
\cnewtheorem{corollary}[theorem]{Corollary}
\cnewtheorem{prop}[theorem]{Proposition}
\cnewtheorem{lemma}[theorem]{Lemma}
\theoremstyle{remark}
\newtheorem{example}{Example}
\newtheorem*{definition}{Definition}
\newtheorem*{notation}{Notation} 
\newtheorem*{remark}{Remark}

\makeatother  %  move after \newtheorem block

\input xypic
\let\xysavmatrix\xymatrix
\def\xymatrix{\disablesubscriptcorrection\xysavmatrix}
\AtBeginDocument{\let\bar\wbar\let\tilde\wtilde\let\hat\what}

\input diagxy
\let\to\rightarrow
\let\barrsquare\square
\let\square\undefined

\makeop{End}

\newcommand{\fqv}{\mathbb{F}_q[V]}
\newcommand{\fqu}{\mathbb{F}_q[u]}

\newcommand{\Euler}{\mbox{\rm\bf E}}
\newcommand{\euler}{\mbox{\rm\bf e}}
\newcommand{\saL}{{\mbox{\rm L}}}
\newcommand{\Mat}{\mbox{\rm Mat}}
 
\newcommand{\Span}{\mbox{\rm Span}} 
\newcommand{\PP}{\mathcal P}
\newcommand{\tr}{{\rm tr}}              	
\newcommand{\ba}{\begin{array}}
\newcommand{\ea}{\end{array}}
\newcommand{\N}{\mathbb N}
\newcommand{\giantP}{{\mathcal P}}
\renewcommand{\P}{{\mathcal P}}
\newcommand{\F}{\mathbb F}
\renewcommand{\O}{{\mathbb O}}
\newcommand{\PF}{{\mathbb P\mathbb F}}
\newcommand{\SO}{{\mathbb S\mathbb O}}
\newcommand{\GL}{\mbox{\rm GL}}
\newcommand{\SL}{\mbox{\rm SL}}
\newcommand{\PSL}{\mbox{\rm PSL}}
\newcommand{\SymMat}{{\rm Mat}^{\rm sym}}
\newcommand{\Sq}{{\mbox{\rm Sq}}}
\def \DD{{\mathbf D}}			
\newcommand{\lra}{\longrightarrow}
\newcommand{\hra}{\hookrightarrow}
\renewcommand{\Im}{\mbox{\rm Im}}

\begin{document}

\begin{abstract} 
The purpose of these notes is to provide an introduction to the Steenrod
algebra in an algebraic manner avoiding any use of cohomology operations.
The Steenrod algebra is presented as a subalgebra of the algebra of
endomorphisms of a functor.  The functor in question assigns to a vector
space over a Galois field the algebra of polynomial functions on that
vector space: the subalgebra of the endomorphisms of this functor that
turns out to be the Steenrod algebra if the ground field is the prime
field, is generated by the homogeneous components of a variant of the
Frobenius map.
\end{abstract}

\maketitle

%%% Document itself goes here.

Beginning with the paper of Adams--Wilkerson \cite{AW}
the Steenrod algebra has played
a significant role in the development of the invariant theory of finite
groups over finite fields (see Smith
\cite{poly,survey}, Neusel--Smith \cite{cyclop} and Neusel \cite{athom} 
and their reference lists). The literature on the Steenrod 
algebra, in particular its
construction, is largely of an algebraic topological nature, making it
difficult for non algebraic topologists to gain insights into how and why 
it
is of significance for invariant theory.
This presents a challenge for those versent in the Steenrod
to explain it to the nonexperts 
in a concise, motivated, and nontechnical\footnote{No
Eilenberg--MacLane spaces, no $\cup_1$ products, etc.}
algebraic manner. More than a decade ago as a visiting professor
at Yale I was confronted with this problem when teaching a course on
invariant theory to an audience consisting
primarily of algbraicists, group theorists, and
number theorists. My strategy to define the Steenrod algebra for this
audience was to
regard the total Steenrod operation as a perturbation of the Frobenius 
map, and
to define the Steenrod algebra as the subalgebra generated by the 
homogeneous
components of this perturbation in the endomorphism algebra of a
carefully chosen functor.

The purpose of these notes is to expand somewhat on that approach and
provide a more complete introduction to the Steenrod
algebra in this manner, ie, 
presented as a subalgebra of the algebra of endomorphisms of a functor.
The functor in question
assigns to a vector space over a Galois field the algebra of
polynomial functions on that vector space: the subalgebra of the
endomorphisms of this functor that turns out to be the Steenrod algebra if
the ground field is the prime field,
is generated by the homogeneous components of a variant of
the Frobenius map.

The material presented here is not new: in fact most of the ideas go
back to the middle of the last century, and are to be found in papers of
H Cartan \cite{HCarone}, \cite{HCar}, J-P Serre \cite{Sefour}, R Thom
\cite{Thom} and Wu Wen-ts\"un \cite{Wu}, with one final key ingredient
being supplied by S Bullet and I Macdonald \cite{BulMac} (see also
T\,P Bisson \cite{bisson}).  My contribution, if there is one, is to
reorganize the presentation of this material so that no algebraic topology
is used, nor is it necessary to assume that the ground field is the prime
field. This way of presenting things appeared \footnote{The emphasis of
\cite[Chapter 10]{poly} is on certain topological applications.} spread
through \cite[Chapters 10 and 11]{poly}. In summary form it also appeared
in \cite{ideals}.  For these notes this material is collected in somewhat
altered form, stripped of its applications to algebraic topology, and
expanded to include the Hopf algebra structure of the Steenrod algebra
due to J\,W Milnor \cite{milnor} for the prime field. For a discussion
of the group of units of the Steenrod algebra regarded as a Hopf algebra
from this point of view see Smith \cite{autos}.

These notes provide a minimal introduction to the use of the Steenrod
algebra in modular invariant theory. The reader is
encouraged to consult the vast literature on the Steenrod algebra.
For orientation in this morass
the reader can do no better than to refer to
%% the notes from Prof. Wood's lectures at the summer school and/or 
the excellent survey article \cite{regwood} and Summer School course notes
\cite{reghits,ingo} by R\,M\,W Wood.

In what follows we adhere to the notations and terminology
of \cite{poly} and \cite{cyclop}.
In particular, if $\F$ is a field and $V = \F^n$ is an $n$--dimensional 
vector
space over $\F$, then $\F[V]$ denotes the graded
algebra of polynomial functions\ on
$V$. This may be regarded as the symmetric algebra on the dual vector 
space
$V^*$ of $V$, where the elements of $V^*$, the linear 
forms have degree 1.
Note carefully 
we ignore the usual topological sign conventions, since graded commutation
rules play no role here. (For a discussion of gradings see eg
\cite[Appendix A Section 1]{cyclop}.)
The correspondence $V \leadsto \F[V]$
defines a contravariant functor from vector spaces over $\F$ to 
graded connected algebras. This functor is at the center of what follows.

\section{The Steenrod algebra}
\label{sec:sec1}
%\sec{steenrod}

We fix once and for all a Galois field
$\F_q$ of characteristic $p$  containing $q = p^\nu$ elements.
Denote by $\F_q[V] [[\xi]]$ the power series ring over $\F_q[V]$
in an additional variable $\xi$, and set $\deg (\xi) = 1-q$.
Define an $\F_q$--algebra homomorphism of degree zero
\[
\giantP (\xi) \co \F_q[V] \lra \F_q[V][[\xi]] ,
\]
by requiring  
\[
\giantP (\xi) (\ell) = \ell + \ell^q\xi \in \F_q[V][[\xi]] , \qua
	\forall  \text{ linear forms }\enspace \ell\in V^*  .
\]
For an arbitrary polynomial $f \in \F_q[V]$, we have after 
separating out homogeneous components,
\footnote{
Let me emphasize here, that we will have no
reason to consider nonhomogeneous polynomials, and implicitly, we are 
always
assuming, unless the contrary is stated, that all algebras are graded, and
if nonnegatively graded, also connected. The algebra
$\F[V][[\xi]]$ is graded, but no longer connected.
}
\begin{equation}
\giantP(\xi)(f) =
\begin{cases}\sum\limits_{i=0}^{\infty} \P^i(f)\xi^i & q \neq 2\\
        	\sum\limits_{i=0}^{\infty} \Sq^i(f)\xi^i& q=2
\end{cases}
\label{ding66}
\end{equation}
This defines $\P^i$, resp. $\Sq^i$, as $\F_q$--linear maps
\[
\P^i , \Sq^i \co \F_q[V]\lra \F_q[V]  .
\]
These maps are functorial in $V$. The operations $\P^i$,
respectively $\Sq^i$, 
%\bold\boldindex{Steenrod squaring
%operations}\bold
are called  \emph{ Steenrod reduced power operations\/}, respectively 
\emph{ Steenrod squaring operations\/}, or collectively, 
\emph{ Steenrod operations\/}. In order to avoid a separate notation for the case 
$q=2$, with the indulgence of topologists, 
\footnote{This is  not the usual
topological convention, which would be to set $\P^i = \Sq^{2i}$. This is 
only
relevant for this algebraic approach when it is necessary to bring in a
Bockstein operation.} we set $\Sq^i = \P^i$ for all $i\in \N_0$.

The sums appearing in \eqref{ding66} are actually finite. In fact
$\giantP(\xi)(f)$ is a {\it polynomial}
in $\xi$ of degree $\deg(f)$ with leading coefficient $f^q$.
This means the Steenrod operations acting on $\F_q[V]$ satisfy 
the \emph{ unstability condition\/}
\[
\P^i(f) = 
\begin{cases}
f^q & i=\deg(f) \\
                  0  & i>\deg(f) \\
\end{cases}\qua \forall   f \in 
\F_q[V].
\]
Note that these conditions express both a triviality condition, viz., 
$\P^i (f)= 0$ for all $i > \deg (f)$, and, a nontriviality condition, viz.,
$\P^{\deg(f)} (f) = f^q$. It is the interplay of these two requirements 
that
seems to endow the unstability condition with the power to yield unexpected
consequences.

Next, observe that the multiplicativity of the operator
$\giantP(\xi)$ leads to  the formulae:
\[
{\P}^k (f' f'') = \sum\limits_{i+j = k} {\P}^i (f') 
{\P}^j (f''), \qua \forall   f', f''\in \F_q[V]  .
\]
These are called the \emph{Cartan formulae\/} for the Steenrod operations.
(NB in field theory, a family of operators satisfying these formulae is
called a \emph{higher order derivation\/}. See, eg Jacobsen
\cite[Chapter 4, Section 9]{jake}.)

As a simple example of how one can compute with these operations consider 
the quadratic form
\[
Q = x^2 + xy + y^2 \in \F_2[x, y]  .
\]
Let us compute how the Steenrod operations $\Sq^i$ act on $Q$ by using
linearity, the Cartan formula, and unstability.
\begin{align*}
\Sq^1 (Q) &= \Sq^1(x^2) + \Sq^1(xy) + \Sq^1(y^2)\\
          &= 2x\Sq^1(x) + \Sq^1(x)\cdot y + x\cdot\Sq^1(y) + 2y\Sq^1(y)\\
          &= 0 + x^2y + xy^2 + 0 = x^2y + xy^2\\
\Sq^2 (Q) &= Q^2 = x^4 + x^2y^2 + y^4\\
\Sq^i (Q) &= 0\text { for } i > 2.
\end{align*}
Since the Steenrod operations are natural with respect to linear
transformations between vector spaces they induce endomorphisms of the 
functor
\[
\F_q[-] \co {\it Vect}_{\F_q} \lra {\it Alg}_{\F_q}
\]
from $\F_q$--vector spaces to commutative graded
$\F_q$--algebras.

They therefore
commute with the action of $\GL(V)$ on $\F_q[V]$.
If $G \hra \GL(n, \F_{q})$ is a faithful representation of a finite group
$G$ then the Steenrod operations restrict to the ring
of invariants $\F_q[V]^{G}$, ie map invariant forms to invariant 
forms.
Hence they can be used to produce new invariants 
from old ones.
This is a new feature of invariant theory over finite fields as opposed to
arbitrary fields (but do see in this connection Glenn \cite{glenn}).
Here is an example to illustrate this. It is based on a result, and the
methods of \cite{segre}.

\begin{example} Let $\F_q$ be the Galois field with 
$q$ elements of odd 
characteristic $p$, and consider the action of the group 
$\SL(2, \F_q)$ on the space of binary quadratic forms over $\F_q$
\
by change of variables. A typical such form is
$Q(x, y) = ax^2 + 2bxy + cy^2$.
\[
\mbox{\rm\bf T}_Q = \left[\ba{cc}
a & b\cr b & c\ea\right]
\]
The space of such forms can be
identified with the vector space $\SymMat_{2, 2}(\F_q)$ of $2 \times 2$
symmetric matrices over $\F_q$.
Under this identification the form $Q$
corresponds to the matrix $\mbox{\rm\bf T}_Q$
above, and the action of $\SL(2, \F_q)$ is
given by $\mbox{\rm\bf T}_Q \mapsto \mbox{\rm\bf S}\mbox{\rm\bf 
T}_Q\mbox{\rm\bf S}^{\tr}$, where
%% $\mbox{\rm\bf T}_Q$ corresponds to the form $Q$, and 
$\mbox{\rm\bf S} \in \SL(2, \F_q)$,
with $\mbox{\rm\bf S}^{\tr}$ the transpose of $\mbox{\rm\bf S}$. 
The element $-\mbox{\rm\bf I} \in \SL(2, \F_q)$ acts trivially. By 
dividing out the
subgroup it generates, we receive a faithful representation
of $\PSL(2, \F_q) = \SL(2, \F_q)\big/\{ \pm\mbox{\rm\bf I}\}$ on the
space of binary quadratic forms. This group has order
$q(q^2 - 1)/2$.

The action of $\PSL(2, \F_q)$ on $\SymMat_{2, 2}(\F_q)$ preserves the
nonsingular quadratic form defined by
$
\det \co \SymMat_{2, 2}(\F_q) \lra \F_q
$
and since there is only one such nonsingular quadratic form in $3$ 
variables
\
over $\F_q$, at least up to isomorphism,
(cf Dickson \cite[pages 169--173]{LED}), we receive an unambiguous faithful
representation $\rho : \PSL(2, \F_q) \hra \O(3, \F_q)$.
Denote by
\[
\left[\ba{cc}
x & y\cr 
y & z\ea\right] \in \SymMat_{2, 2}(\F_q)^*
\]
a generic linear form on the dual space of the $2 \times 2$ symmetric
matrices over
$\F_q$. Per definition the quadratic form
\[
\det = xz - y^2 \in \F_q[\SymMat_{2, 2}(\F_q)] = \F_q[x, y, z]
\]
is $\O(3, \F_q)$--invariant. If we apply the first Steenrod
operation to this form we receive the new invariant form of degree $q+1$,
viz.,
\[
\P^1 (\det) = x^qz + xz^q - 2y^{q+1} \in \F_q[x, y, z]^{\O(3, \F_q)}  .
\]

The full ring of invariants of the orthogonal group $\O(3, \F_q)$ is
known (see eg Cohen \cite{sdcohen} or \cite{segre}).
To wit
\[
\F_q[x, y, z]^{\O(3, \F_q)} \cong  \F_q[\det, \P^1(\det), \Euler_{\det}]
  .
\]
Here $\Euler_{\det}$ is the Euler class (see eg Smith--Strong \cite{orbits} 
or Neusel--Smith \cite[Chapter 4]{cyclop})
associated to the
configuration of linear 
forms defining the set of external lines to the projective variety
${\mathfrak X}_{\det}$ in the projective plane
$\PF_q(2)$ over $\F_q$
defined by the vanishing of the
determinant
\footnote{
The projective plane of $\F_q$ is defined by
$\PF_q(2) = \big(\F_q^3 \setminus \lbrace 0\rbrace \big)^{\F^\times}$
where $\F^{\times}$ acts via scalar multiplication on the vectors of $\F_q^3$.
In this discussion we are identifying $\F_q^3$ with
$\Mat_{2, 2}^{\rm sym} (\F_q)$, so
this is the same as the set of lines through the origin in
$\Mat_{2, 2}^{\rm sym}$. The pre-Euler class $\euler_{\det}$
may be taken to be the product of a set of linear forms
$\big\lbrace \ell_{\mbox{\scriptsize \saL}}\big\rbrace$, 
indexed by the 
$\binom{q}{2}$ 
external lines 
$\lbrace{\saL}\rbrace$ to 
${\mathfrak X}_{\det}$,
and satisfying
$\ker(\ell_{\mbox{\scriptsize \saL}}) = {\saL}$. The Euler class 
$\Euler_{\det}$ is its
square.
}
(see Hirschfeld \cite[Section 8.2]{hirschfeld} and \cite{segre}).
The form $\Euler_{\det}$ has degree $q(q - 1)$.  The three forms
$\det, \P^1(\det), \Euler_Q \in \F_q[x, y, z]^{\O(3, \F_q)}$
are a system of parameters \cite{segre}. Since
the product of their degrees is $|\O(3, \F_q)|$ it follows from 
\cite[Proposition 5.5.5]{poly} that $\F_q[x, y, z]^{\O(3, \F_q)}$ must be a 
polynomial algebra as stated.

The pre-Euler class $\euler_{\det}$ of the set of external projective 
lines to
${\mathfrak X}_{\det}$ is an $\O(3, \F_q)$ $\det$--relative invariant, so
is $\SO(3, \F_q)$--invariant. It has degree $\binom{q}{2}$, and together 
with
the forms $\det$ and $\P^1(\det)$ it forms a system of parameters for
$\F_q[x, y, z]^{\SO(3, \F_q)}$, so again we may apply 
\cite[Proposition 5.5.5]{poly}
and conclude that $\F_q[x, y, z]^{\SO(3, \F_q)}$ is a polynomial algebra, 
viz.,
$\F_q[x, y, z]^{\SO(3, \F_q)} = \F_q[\det, \P^1(\det), \euler_{\det}]$.

Finally, $\PSL(2, \F_q)$ is the commutator subgroup of, and has index $2$
in,
the special orthogonal group
$\SO(3, \F_q)$,
so by a Proposition in Smith \cite{twobytwo}
the ring of invariants of $\PSL(2, \F_q)$ acting on the space of binary
quadratic forms is a hypersurface. It has generators
$\det, \P^1(\det), \euler_{\det}$ and a certain form $\omega$ which 
satisfies a
monic quadratic equation over the subalgebra generated by the first three.
A choice for $\omega$ is the pre-Euler class of the configuration of 
external
projective lines to the variety ${\mathfrak X}_{\det} \subset \PF(3)$.
\end{example}

The Steenrod operations can be collected together to form an algebra,
in fact a Hopf algebra (see \fullref{sec:sec4}), 
over the Galois field $\F_q$.

\begin{definition}
The \emph{Steenrod algebra\/} $\PP^*(\F_q)$ 
is the $\F_q$--subalgebra of the endomorphism algebra of the functor
$\F_q[-]$, generated by  $\P^0= 1, \P^1, \P^2 ,\dots$.
\end{definition}

\begin{notation} In most situations, such as here,
the ground field $\F_q$ is fixed at the outset, and we
therefore abbreviate $\PP^*(\F_q)$ to $\PP^*$.
\end{notation}

The next sections develop the basic algebraic structure
of the Steenrod algebra.

\section{The Adem--Wu relations}
\label{sec:sec2}
%\sec{relations}

The Steenrod algebra is by no means freely generated by the Steenrod
reduced powers.
For example, when $p=2$ it is easy to check that
$\Sq^1\Sq^1 = 0$ by verifying this is the case for monomials
$z^E = z_1^{e_1},\ldots , z_n^{e_n}$.  To do so one needs the formula,
valid for any linear form,
$\Sq^1 (z^k) = k z^{k+1}$, 
which follows by induction immediately
from the Cartan formula.\footnote{In fact every
element in the Steenrod algebra is nilpotent, but the index of nilpotence 
is
known only in a few cases, see eg Monks \cite{monksone,monkstwo},
Walker--Wood \cite{WWone,WWtwo} and Wood
\cite{regwood} for a resum\'e of what is known.}

Traditionally, relations between the Steenrod operations are expressed as
commutation rules for ${\cal P}^{i}{\cal P}^{j}$, respectively
$\Sq^{i}\Sq^{j}$. These commutation relations are called
\emph{ Adem--Wu relations\/}.
In the case of the prime field $\F_p$\
they were originally conjectured
by Wu based on his study of the mod p cohomology of Grassmann
manifolds \cite{Wu} and proved by J Adem in \cite{adem}, H Cartan in
\cite{HCarone}, and for $p=2$ by J\,P Serre in \cite{Sefour}.
These relations are usually written as follows:
\[
{\P}^i {\P}^j =  \sum_{k=0}^{[i/q]} (-1)^{(i-qk)} \cdot
	{\binom{(q-1)(j-k)-1}{i-qk}} {\P}^{i+j-k}
{\P}^k\quad \forall   i, j \geq 0, i<qj.
\]
Note for any Galois field $\F_q$ the coefficients are still elements in 
the 
prime subfield $\F_p$ of $\F_q$. 

The proof of these relations is greatly simplified by the
%\bold
\emph{Bullett--Macdonald identity\/}, which provides us with a 
well-wrapped 
description of the relations among the Steenrod operations,
Bullett--Macdonald \cite{BulMac}.
To describe this identity, as in \cite{BulMac},
extend $\giantP(\xi)$ to a ring homomorphism
$\giantP(\xi) \co \F[V][\eta] \lra \F[V][\eta][\xi]$
by setting $\giantP(\xi))(\eta)=\eta$.
Next, set $u = (1-t)^{q-1} = 1 + t + \cdots + t^{q-1}$ and  $s = tu$. 
Then
the Bullett--Macdonald identity is
\[
\giantP(s)\circ \giantP(1) = \giantP(u)\circ\giantP(t^q).
\] 
Since $\giantP(\xi)$ is additive and multiplicative, it is enough
to check this equation for the basis elements of $V^*$, which is indeed a 
short calculation. 
Rumor says Macdonald, like most of us, could not remember the coefficients
that appear in the Adem relations, so devised this identity so that he 
could
derive them on the spot when J\,F Adams came to talk with him.

\begin{remark} For $p=2$, T\,P Bisson has pointed out 
(see Bisson--Joyal \cite{bissonjoyal})
that the
Bullett--Macdonald identity may be viewed as a commutation rule, viz.,
$\giantP(\xi)\giantP(\eta) = \giantP(\eta)\giantP(\xi)$.
For a general Galois $\F_q$, one needs to demand 
$\GL(2, \F_q)$--invariance
of $\giantP(\zeta)$, where
$\zeta \in \Span_{\F_q}\{ \xi, \eta\}$.
\end{remark}

To derive the Adem--Wu relations we
provide details for the residue computation\footnote{The following 
discussion
is based on conversations with E\,H Brown Jr. The author is also grateful 
to J Hartmann for correcting some errors in his version of the
computation. I do hope for once the indices are close to being correct.}
%% 22.10.2003 : optional hyphenation added
sketched in \cite{BulMac}.
First of all, direct calculation gives:
\begin{align*}
\giantP(s)\giantP(1) &= \sum_{a,\ k} s^a \P^a \P^k\\
\giantP(u)\giantP(t^q) &=\sum_{a,\ b,\ j} u^{a+b-j} t^{qj} \P^{a+b-j} \P^j,
\end{align*}
which the Bullett--Macdonald identity
says are equal. Recall from complex analysis that
\[
\frac{1}{2\pi \imath} \oint_\gamma z^m dz =
\begin{cases} 1 & m = -1\\ 
0 & \text{otherwise},
\end{cases}
\]
where $\gamma$ is a small circle around $0 \in \C$. Therefore we obtain
\begin{align*}
\sum_k \P^a \P^k &= \frac{1}{2\pi\imath} 
	\oint_\gamma \frac{\giantP(s) \giantP(1)}{s^{a+1}}ds\\
&= \frac{1}{2\pi\imath} \oint_\gamma \frac{\giantP(u) 
\giantP(t^q)}{s^{a+1}} ds\\
&= \frac{1}{2\pi\imath}
\sum_{a,\ b,\ j} \oint_\gamma \frac{u^{a+b-j}t^{qj}}{s^{a+1}} ds
	\P^{a+b-j} \P^j.
\end{align*}
The formula $s = t(1 - t)^{q-1}$ gives
$ds = (1 - t)^{q-2} (1 - qt) dt$, so
%% putting $c = a+b$ and
substituting gives
\begin{align*}
\frac {u^{a+b-j} t^{qj}} {s^{a+1}} ds
&= \frac{ (1-t)^{(q-1)(a+b-j)} t^{qj} (1-t)^{q-2} (1-qt) }
	{\left[t (1-t)^{q-1}\right]^{a + 1} } dt\cr
&= (1-t)^{(b-j-1)(q-1) + (q-2)} t^{qj-a-1} (1 - qt) dt\cr
&= (1-t)^{((b-j)(q-1) -1)} t^{qj-a-1} (1 - qt) dt\cr
&= \left[ \sum_{k} (-1)^k { \binom{(b-j)(q-1) -1}{k} } t^k \right]
		t^{qj -a-1} (1 - qt) dt\cr
&= \sum_{k} (-1)^k { \binom{(b-j)(q-1) - 1} {k} }
	\left[ t^{k + qj - a -1} - q t^{k + qj - a} \right] dt.
\end{align*}
Therefore
\begin{multline*}
\P^a \P^b = \sum_j \left[ \frac{1}{2\pi\imath} \oint_\gamma
	\frac{u^{a{+}b{-}j} t^{qj}} {s^{a+1}} ds\right] \P^{a{+}b{-}j} \P^j \\
 = \sum_j \frac {1}{2\pi\imath}\oint 
	\sum_{k} (-1)^k { \binom{(b{-}j)(q{-}1){-}1}{k} }
\left[ t^{k{+}qj{-}a{-}1}{-}qt^{k{+}qj{-}a} \right] dt \P^{a{+}b{-}j} \P^j.
\end{multline*}
Only the terms where
\begin{eqnarray*}
k + qj -a -1  &=& -1 \quad (k = a - qj)\cr
k + qj -a  &=& -1 \quad (k = a - qj - 1),
\end{eqnarray*}
contribute anything to the sum, so
\begin{align*}
\P^a \P^b & =
\sum_j \biggl[ (-1)^{(a-qj)}\cdot{\binom {(b-j)(q-1) - 1} {a - qj} } \\
& \qua +
	(-1)^{a - qj - 1} q { \binom{(b-j)(q-1) -1}  {a - qj - 1} } 
\biggr]
		\P^{a+b-j} \P^{j},
\end{align*}
and since
\[
{ \binom{(b-j)(q-1) - 1} {a - qj} } -
	q { \binom{(b-j)(q-1) -1} {a - qj - 1} }
	\equiv { \binom{(b-j)(q-1) - 1}  {a - qj} } \bmod p,
\]
we conclude
\[
\P^a \P^b = \sum_j (-1)^{(a-qj)}\cdot{\binom{(b-j)(q-1) - 1} {a - qj} }
	\P^{a+b-j} \P^{j},
\]
which are the Adem--Wu relations.

Thus there is a surjective map from the free associative algebra 
${\cal B}^{*}$ with 1
generated by the Steenrod operations modulo the ideal generated by the
Adem--Wu relations, %%which we denote by ${\cal B}^{*}$,
\[
\P^a \P^b - \sum_j (-1)^{(a-qj)}\cdot{\binom{(b-j)(q-1) - 1} {a - qj} }
	\P^{a+b-j} \P^{j}\quad a, b \in \N\ \hbox{\rm and}\ a < qb, 
\]
onto the Steenrod algebra.
In fact, this map, ${\cal B}^{*} \lra {\cal P}^{*}$
is an isomorphism, so
the Adem--Wu relations are a complete set of defining relations for the
Steenrod algebra. The proof of this, and some of its consequences,
is the subject of the next section.

\section{The basis of admissible monomials}
\label{sec:sec3}
%\sec{basis}

In this section we show that the relations
between Steenrod operations that are universally valid all follow from
the Adem--Wu relations. To do so
we extend some theorems of
of H Cartan \cite{HCarone}, 
J-P Serre \cite{Sefour} and Wu Wen-ts\"un \cite{Wu}
from the case of the prime field to
arbitrary Galois fields. Their proofs have been rearranged so that
no direct use is made of algebraic topology.

An \emph{ index sequence\/} is a sequence 
$I = (i_1, i_2, \ldots, i_k,\ldots)$ of nonnegative
integers, almost all of which are zero. If $I$ is an index sequence
we denote by
$\P^I \in \PP^*$ the monomial 
$\P^{i_1}\cdot\P^{i_2}\cdots\P^{i_k}\cdots$
in the Steenrod operations $\P^i$, with the convention that
trailing $1$s are ignored. The degree of the element $\P^I$ is
$(q-1)(i_1 + i_2 + \cdots + i_k +\cdots)$.
These iterations of Steenrod operations are called \emph{ basic monomials\/}.
An index sequence $I$ is called
admissible if $i_s \geq q i_{s+1}$ for $s = 1,\dots$.
We call $k$ the \emph{ length\/}
of $I$ if $i_k \ne 0$ but $i_s = 0$ for $s > k$.
Write $\ell (I)$ for the length of $I$.
It is often convenient to treat an index sequence as a finite sequence of
nonnegative integers by truncating it to $\ell(I)$ entries.

A basic monomial is defined to be \emph{ admissible\/}
if the corresponding index sequence is admissible. The strategy of
H Cartan and J-P Serre to prove
that the Adem--Wu relations
are a complete set of defining relations for the Steenrod algebra of the
prime field is to prove
that the admissible monomials are an $\F_p$ basis for $\PP^*$. We follow
the same strategy for an arbitrary Galois field.

Recall that ${\cal B}^*$ denotes the free, graded, associative algebra 
generated
by the symbols $\P^k$ modulo the ideal spanned by the Adem--Wu relations
in those symbols.
We have a surjective map ${\cal B}^* \lra \PP^*$, and so with
his notation our goal is to prove:

\begin{theorem}
The admissible monomials span 
${\cal B}^{*}$ as an
$\F_{q}$--vector space.
The images of the admissible monomials in the Steenrod algebra
are linearly independent.
\label{thm:thm3.1}
\end{theorem}

\begin{proof}
We begin by showing that the 
admissible monomials span ${\cal B}^{*}$.

For a sequence $I = (i_{1},i_{2},\ldots , i_{k}) $,
the \emph{ moment\/} of $I$, denoted by $m(I)$, is defined 
by $m(I) = \sum_{s=1}^{k}s\cdot i_{s}$. We first show that an inadmissible
monomial is a sum of monomials of smaller moment.
Granted this it follows by induction over the moment that the admissible
monomials span ${\cal B}^{*}$.

Suppose that $\P^{I}$ is an inadmissible monomial. Then there is a
smallest $s$ such
that $i_{s} < qi_{s+1}$, ie,
\[
\P^{I} = {\cal Q}'\P^{i_{s}}\P^{i_{s+1}}{\cal Q}'',
\]
where ${\cal Q}'$, ${\cal Q}''$ are basic monomials, and ${\cal Q}'$ is
admissible. It is therefore
possible to apply an Adem--Wu relation to $\P^{I}$ to obtain
\[
\P^{I} = 
\sum_j a_{j} {\cal Q}'\P^{i_{s} + i_{s+1} - j}\P^{j}{\cal Q}'',
\]
for certain coefficients $a_{j} \in \F_{p}$. The terms on the right hand
side all have smaller moment than $\P^{I}$ and so, by induction on $s$,
we may express $\P^{I}$ as a sum of admissible monomials.
(NB The admissible
monomials are {\it reduced} in the sense that no Adem--Wu relation can be 
applied to them.)

We next show that the admissible monomials are linearly independent as 
elements
of the Steenrod algebra $\PP^*$. This we do by adapting
an argument of J-P Serre \cite{Sefour} and H Cartan \cite{HCarone}
which makes use of a formula of Wu Wen-ts\"un.

Let $e_{n} = x_{1}x_{2}\cdots x_{n} \in \F_{q}[x_{1},\ldots ,  x_{n}]$
be the $n$th elementary symmetric function.
Then,
\begin{eqnarray*}
\giantP(\xi)(e_{n}) &=& \giantP(\xi)(\prod_{i=1}^{n} x_{i})
	= \prod_{i=1}^{n} \giantP(\xi)(x_{i})\cr
		&=& \prod_{i=1}^{n} (x_{i} + x_{i}^{q}\xi)
	= \prod_{i=1}^{n} x_{i}\cdot\prod_{i=1}^{n}(1 + x_{i}^{q-1}\xi)\cr
	&=& e_n(x_{1},\ldots , x_{n})\cdot
\left(\sum_{i=1}^{n} e_{i}(x_{1}^{q-1},\ldots ,  x_{n}^{q-1})\xi^i\right)
 ,
\end{eqnarray*}
where $e_i(x_1,\ldots , x_n)$ denotes the $i$th elementary symmetric
polynomial in $x_1,\ldots , x_n$.
So we have obtained a formula
of Wu Wen-ts\"un:
\[ \P^{i}(e_{n}) = e_{n}\cdot e_{i}(x_{1}^{q-1},\ldots , x_{n}^{q-1}) \]
We claim that the monomials
\[ \{ \P^{I} | \P^I \text{ admissible and } \deg(\P^{I}) \leq 2n \} \]
are linearly independent in $\F_{q}[x_{1},\ldots ,  x_{n}]$. To see this
note that in case $\ell (I) \leq n$, each entry in $I$ is at
most $n$ (so the following formula makes sense), and
\[
\P^{I}(e_{n}) = e_{n}\cdot
	\prod_{j = 1}^{s} e_{i_{j}}(x_{1}^{q-1},\ldots ,  x_{n}^{q-1}) + 
\cdots
\]
where $I = (i_{1},\ldots , i_{s})$,
$\P^{I} = \P^{i_{1}}\cdots \P^{i_{s}}$ and
the remaining terms are lower in the lexicographic ordering on monomials.
So $e_{n}\cdot \prod_{j = 1}^{s} e_{i_{j}}(x_{1}^{q-1},\ldots ,  
x_{n}^{q-1})$
is the largest monomial in $\P^{I}(e_{n})$ in the
lexicographic order. Thus
\[
\{\P^{I}(e_n)\ |\ \P^I\ {\rm admissible\ and} \ 
	\deg(\P^{I}) \leq 2n\} ,
\]
have distinct largest monomials, so are linearly independent.

By letting $ n \lra \infty $ we obtain the assertion, completing the proof.
\end{proof}

Thus the
Steenrod algebra may be regarded (this is one traditional definition)
as the graded
free associative algebra
with 1 generated by the $\Sq^{i}$ respectively $\P^{i}$ modulo the ideal
generated by the Adem--Wu relations. This means we have proven:

\begin{theorem} The Steenrod algebra $\PP^*$ is 
the free
associative $\F_q$--algebra  generated by the reduced power operations
$\P^0, \P^1, \P^2, \dots$ modulo the Adem--Wu relations.
\label{thm:thm3.2}
\end{theorem}

\begin{corollary}
The admissible monomials are 
an $\F_q$--basis for the Steenrod algebra $\PP^*$.
\label{cor:cor3.3}
\end{corollary}

Since the coefficients of the Adem--Wu relations lie in the prime field
$\F_p$, the operations $\P^{p^i}$ for $i \geq 0$
are indecomposables in $\PP^*$. In particular, over the Galois field
$\F_q$, the Steenrod algebra
$\PP^*$ is \emph{ not\/} generated by the operations $\P^{q^i}$ for 
$i \geq 0$:
one needs all $\P^{p^i}$ for $i \geq 0$.
This will become even clearer after we have developed the Hopf algebra
structure of $\PP^*$ in the next section.

\begin{example} Consider the polynomial algebra 
$\F_2[Q, T]$
over the field with $2$ elements, where the indeterminate $Q$ has degree 
$2$
and $T$ has degree $3$. If the Steenrod algebra were to act unstably on 
this
algebra then the unstability condition would determine $\Sq^i (Q)$ and
$\Sq^j(T)$ apart from $i=1$ and $j = 1$ and $2$. If we specify these as
follows
\[
\Sq^1 (Q) = T,\qua \Sq^1 (T) = 0,\qua \Sq^2(T) = QT ,
\]
and demand that the Cartan formula hold,
then using these formulae we can compute $\Sq^k$ on any monomial, and hence
by linearity, on any polynomial in $Q$ and $T$. For example
\[
\Sq^1(QT) = \Sq^1(Q)\cdot T + Q\cdot\Sq^1(T) = T^2 + 0 = T^2 ,
\] 
and so on. 
Note that since
$\Sq^1\cdot\Sq^1 = 0$ is an Adem--Wu relation, $\Sq^1(T) = 0$ is forced
from $\Sq^1(Q) = T$. To verify the unstability conditions, suppose that
\[
\Sq^a\Sq^b = \sum_{c=0}^{[\frac{a}{2}]} {\binom{b-1-c}{a-2c}}
	\Sq^{a+b-c}\Sq^{c} ,\qua 0 < a < 2b ,
\]
is an Adem--Wu relation. We need to show that
\[
\left(\Sq^a\Sq^b - \sum_{c=0}^{[\frac{a}{2}]} {\binom{b-1-c}{a-2c}}
	\Sq^{a+b-c}\Sq^{c}\right)\big(Q^iT^j\big) = 0
\]
for all $i, j \in \N_0$.
By a simple argument using the Cartan formulae, see eg
\cite[Lemma 4.1]{steenrod}, it is enough to verify that
these hold for the generators $Q$ and $T$ and this is routine.
It is a bit
more elegant to identify $Q$ with $x^2 +xy + y^2$
and $T$ with $x^2y + xy^2 \in \F_2[x, y]$. The action of the Steenrod
operations on $Q$ and $T$ then coincides with the restriction of the
action from $\F_2[x, y]$. This way, it is then clear
that $\F_2[Q, T]$ is
an unstable algebra over the Steenrod algebra, because,
\begin{enumerate}
\item[(1)] with some topological background we recognize this as just 
$H^*(B\SO(3); \F_2)$ , or,
\item[(2)] with some invariant theoretic background we recognize this
as the Dickson algebra
$\DD(2) = \F_2[x, y]^{\mbox{\scriptsize{\GL}}(2, \F_2)}$.
\end{enumerate}
\end{example}

\section{The Hopf algebra structure of the Steenrod algebra}
\label{sec:sec4}

%\sec{hopf}
Our goal in this section is to complete the traditional picture of the 
Steenrod
algebra by proving that
$\PP^*(\F_q)$ is a Hopf algebra
\footnote{One quick way to do this is to write
down as comultiplication map
\[
\nabla(\P^k) = \sum_{i + j = k} \P^i \otimes \P^j ,\quad k = 1,
2, \dots ,
\]
and verify that it is compatible with the Bullett--Macdonald identity, and 
hence
also with the Adem--Wu relations.}
and extending Milnor's Hopf algebra
\cite{milnor} structure theorems from the prime field 
$\F_p$ to an arbitrary
Galois field. It should be emphasized that this requires no new ideas, 
only a
careful reorganization of Milnor's proofs, so as to
avoid reference to algebraic topology and
cohomology operations, and, where appropriate
carefully replacing $p$ by $q$.

%% prove it is a Hopf algebra

\begin{prop}
Let $p$ be a prime integer, 
$q=p^\nu$ a power of $p$, and
$\F_q$ the Galois field with $q$ elements. Then the Steenrod algebra of 
$\F_q$
is a cocommutative Hopf algebra over $\F_q$ with respect to the coproduct
\[
\nabla\co \PP^* \lra \PP^* \otimes \PP^*,
\]
defined by the formulae
\[
\nabla(\P^k) = \sum_{i + j = k} \P^i \otimes \P^j,\quad k = 1,
2, \dots .
\]
\end{prop}

\begin{proof} 
Consider the functor
$V \leadsto \F_q[V]\otimes \F_q[V]$
that assigns to a finite dimensional
vector space $V$ over $\F_q$ the commutative graded
algebra $\F_q[V]\otimes\F_q[V]$ over $\F_q$.
There is a natural map of algebras
\[
\PP^* \otimes \PP^* \lra {\rm End}(V \leadsto \F_q[V]\otimes \F_q[V]),
\]
given by the tensor product of endomorphisms.
Since there is an isomorphism 
$$\F_q[V]\otimes\F_q[V] \cong  \F_q[V\oplus V],$$
that is natural in $V$, the functor
${\rm End}(V \leadsto \F_q[V]\otimes \F_q[V])$
is a subfunctor of the functor
${\rm End}(V \leadsto \F_q[V])$
that assigns to a finite dimensional vector space 
$V$ over $\F_q$ the polynomial algebra $\F_q[V]$. Hence restriction
defines a map of algebras
\[
\PP^* \lra  {\rm End}(V \leadsto \F_q[V]\otimes \F_q[V]),
\]
and we obtain a diagram of algebra homomorphisms
$$\bfig \dtriangle/<--`->`->/<1200,500>[\PP^*\otimes\PP^*`\PP^*`
  \End(V\leadsto\fqv\otimes\fqv);`\tau`\rho] \efig$$
What we need to show is that $\Im(\rho) \subseteq \Im(\tau)$,
for since $\tau$ is monic $\nabla = \tau^{-1}\rho$ would define the desired
coproduct. Since the reduced power operations
$\P^k$ for $k = 1, 2,\ldots , $ generate $\PP^*$
it is enough to check that $\rho (\P^k) \in \Im(\tau)$
for $k = 1, 2,\ldots ,$. But this is immediate\
from the Cartan formula. Since $\nabla$ is a map of algebras the Hopf
condition is satisfied, so $\PP^*$ is a Hopf algebra.
\end{proof}

If $J$ is an admissible index sequence then
\[
e(J) = \sum_{s=1}^{\infty} (j_s - q j_{s+1}),
\]
is called the \emph{ excess\/} of $J$.
For example, the sequences
\[
M_k = (q^{k-1},\ldots , q, 1),\qua k = 1, 2,\dots,
\]
are all the admissible sequences of excess zero.
Note that 
\[
\deg(\P^{M_k}) = \sum_{j=1}^k q^{k-j} (q - 1) = q^k - 1,
\text{ for } k = 1, 2,\dots  .
\]

Recall by \fullref{cor:cor3.3} that the admissible monomials are an 
$\F_q$--vector space basis for $\PP^*$.

Let $\PP_*(\F_q)$ denote the Hopf algebra dual to the Steenrod algebra
$\PP^*(\F_q)$.\
We define $\xi_k \in \PP_*(\F_q)$ to be dual to the monomial
$\P^{M_k} = \P^{q^{k-1}}\cdots \P^q\cdot\P^1$
with respect to the basis of admissible monomials for $\PP^*$. 
This means that we have:
\[
\langle  \P^J \mid \xi_k \rangle = 
\begin{cases}
		1 & J = M_k\\
		0 & \text{ otherwise},
\end{cases}
\]
where we have written $\langle  \P \mid \xi \rangle$ for the value of an 
element
$\P \in \PP^*(\F_q)$ on an element $\xi \in \PP_*(\F_q)$.
Note that $\deg (\xi_k) = q^k - 1$ for $k = 1,\ldots$.

If $I = (i_1, i_2,\ldots , i_k,\dots)$ is an index sequence
we call $\ell$ the \emph{ length\/} of $I$,
denoted by $\ell(I)$, if $i_k = 0$ for $k > \ell$, but
$i_\ell \neq 0$. We associate to an index sequence
$I = (i_1, i_2,\ldots , i_k,\dots)$ the element
$\xi^I = \xi_{1}^{i_1}\cdot\xi_{2}^{i_2}\cdots\xi_\ell^{i_\ell} \in
\PP_*(\F_q)$, where $\ell = \ell(I)$. Note that
\[
\deg (\xi^I) = \sum_{s=1}^{\ell(I)} i_s (q^s - 1)  .
\]
To an index sequence
$I = (i_1, i_2,\ldots , i_k,\dots)$ we also associate an admissible
sequence $J(I) = (j_1, j_2,\ldots , j_k,\dots)$ defined by
\begin{equation}
j_1 = \sum_{s=1}^\infty i_s q^{s-1},\qua
j_2 = \sum_{s=2}^\infty i_s q^{s-2},\ldots,\qua
j_k = \sum_{s=k}^\infty i_s q^{s-k},\dots  .
\label{eqn:ding67}
\end{equation}
It is easy to verify that as $I$ runs over all index sequences that $J(I)$ 
runs
over all admissible sequences. Finally, note that
$\deg (\P^{J(I)}) = \deg (\xi^I)$ for any index sequence $I$.

The crucial observation used by Milnor to
prove the structure theorem of $\PP_*(\F_q)$
is that the pairing of the admissible monomial basis for $\PP^*(\F_q)$
against the  monomials in the $\xi_k$ is upper triangular.
To formulate this precisely we order the index sequences lexicographically 
from the right so for example
$(1, 2, 0,\dots) \prec (0, 0, 1,\dots)$ .

\begin{lemma}[J\,W Milnor]
With the preceding notations we have that the inner product matrix
$\langle  \P^{J(I)} \mid \xi^K \rangle$
is upper triangular with $1$s on the diagonal, ie,
\[
\langle  \P^{J(I)} \mid \xi^K \rangle = 
\begin{cases}
1 & I = K \\
0 & I < K.
\end{cases}
\]
\label{lem:lem4.2}
\end{lemma}

\begin{proof}
Let the length of $K$ be $\ell$ and define
$K' = (k_1, k_2,\ldots , k_{\ell - 1})$, so
\[
\xi^K = \xi^{K'}\cdot \xi_\ell \in \PP_*(\F_q)  .
\]
If $\nabla$ denotes the coproduct in $\PP^*(\F_q)$, then
we have the formula
\begin{equation}
\langle  \P^{J(I)} \mid \xi^K \rangle =
	\langle  \P^{J(I)} \mid \xi^{K'}\cdot \xi_\ell \rangle =
\langle  \nabla(\P^{J(I)}) \mid \xi^{K'}\otimes \xi_\ell \rangle
\label{eqn:ding64}
\end{equation}
If $J(I) = (j_1, j_2,\ldots , j_k,\dots)$ then one easily checks
that
\[
\nabla(\P^{J(I)}) = \sum_{J' + J'' = J(I)} \P^{J'} \otimes \P^{J''}  .
\]
Substituting this into \eqref{eqn:ding64} gives
\begin{equation}
\langle  \P^{J(I)} \mid \xi^K \rangle =
	\sum_{J' + J'' = J(I)} \langle  \P^{J'} \mid \xi^{K'} \rangle
	\cdot
	\langle  \P^{J''} \mid \xi_\ell \rangle  .
\label{eqn:ding72}
\end{equation}
By the definition of $\xi_\ell$ we have
\[
\langle  \P^{J''} \mid \xi_\ell \rangle = 
\begin{cases}
1 & J'' = M_\ell \\
0 & \text{otherwise}.
\end{cases}
\]
If $J'' = M_\ell$ then unravelling the definitions shows that
$J' = J(I')$, for a suitable $I'$,
so if $K$ and $I$ have the same length $\ell$, we have shown
\[
\langle  \P^{J(I)} \mid \xi^K \rangle =
\langle  \P^{J(I')} \mid \xi^{K'} \rangle ,
\]
and hence it follows from induction over the degree that
\[
\langle  \P^{J(I)} \mid \xi^K \rangle = 
\begin{cases}
1 & I = K\\
0 & I < K.
\end{cases}
\]
If, on the other hand, $\ell (I) < \ell$ then all the terms
\[
	\langle  \P^{J''} \mid \xi_\ell \rangle
\]
in the sum \eqref{eqn:ding72} are zero and hence that
$\langle  \P^{J(I)} \mid \xi^K \rangle = 0$ as required.
\end{proof}

\begin{theorem}
Let $p$ be a prime integer, 
$q=p^\nu$ a power of $p$, and
$\F_q$ the Galois field with $q$ elements. Let $\PP_*(\F_q)$ denote the 
dual
Hopf algebra to the  Steenrod algebra of the Galois field $\F_q$. Then,
as an algebra
\[
\PP_* \cong  \F_q[\xi_1,\ldots , \xi_k,\dots] ,
\]
where $\deg (\xi_k) = q^k - 1$ for $k \in \N$. The coproduct is given by 
the
formula
\[
\nabla_* (\xi_k) = \sum_{i + j = k} \xi_i^{q^j} \otimes \xi_j,\quad
k = 1, 2,\dots  .
\]
\end{theorem}

\begin{proof}
By Milnor's Lemma (\fullref{lem:lem4.2})  the monomials
$\big\lbrace\xi^I\big\rbrace$ where $I$ ranges over all index sequences 
are linearly
independent in $\PP_*(\F_q)$. Hence $\F_q[\xi_1,\ldots , \xi_k,\dots]
\subseteq \PP_*(\F_q)$.
%% But 
% 22/10.2003 "But" replaced by "The algebras" to improve line break
The algebras
$\PP_*(\F_q)$ and
$\F_q[\xi_1,\ldots , \xi_k,\dots]$ have the same Poincar\'e series, since
\
$\deg (\P^{J(I)}) = \deg (\xi^I)$ for all index sequences $I$, and the
admissible monomials $\P^{J(I)}$ are an $\F_q$--vector space basis for
$\PP^*(\F_q)$. So
$\F_q[\xi_1,\ldots , \xi_k,\dots] = \PP_*(\F_q)$, and it remains to
verify the formula for the coproduct.

To this end we use the test algebra $\F_q[u]$, the polynomial algebra on 
one
generator, as in \cite{milnor}. Note that for admissible sequences we have
\begin{equation}
\P^J (u) = \begin{cases}
u^{q^k} & J = M_k \\
0 & \text{otherwise}.
\end{cases}
\label{eqn:ding81}
\end{equation}
Define the map
\[
\lambda^*\co \F_q[u] \lra \F_q[u] \otimes \PP_*,
\]
by the formula
\[
\lambda^* (u^i) = \sum \P^{J(I)} (u^i) \otimes \xi^I,
\]
where the sum is over all index sequences $I$. Note that in any given 
degree the sum is finite and that $\lambda^*$ is a map of algebras. 
Moreover
\[
(\lambda^*\otimes 1)\lambda^* (u) = (1 \otimes \nabla_*)\lambda^*(u) ,
\]
ie the following diagram
\begin{equation}
\bfig
\barrsquare/<-`<-`<-`<-/<1500,500>[
  \fqu\otimes\PP_*(\F_q)\otimes\PP_*(\F_q)`
  \fqu\otimes\PP_*` \fqu\otimes\PP_*` \fqu;
  1\otimes\nabla_*`\lambda^*\otimes 1`\lambda^*`\lambda^*]
  \efig
%\ba{ccc}
%F_q[u]\otimes\PP_*(\F_q)\otimes\PP_*(\F_q) & \vlpf{10}{1\otimes\nabla_*}
%	& \F_q[u] \otimes \PP_*\\
%\opf{\lambda^*\otimes 1}& &\opf{\lambda^*}\\
%	\F_q[u] \otimes \PP_*&
%\hbox to 0pt{\kern-4pc$\vlpf{20}{\lambda^*}$\kern4pc\hss}
%& \F_q[u]
%\ea
%\xymatrix{
%F_q[u]\otimes\PP_*(\F_q)\otimes\PP_*(\F_q) 
%& &
%\F_q[u] \otimes \PP_*  \ar[ll]_{1\otimes\nabla_*}
%\\
%\F_q[u] \otimes \PP_* \ar[u]_{\lambda^*\otimes 1} & &
%\F_q[u] \ar[u]_{\lambda^*} \ar[ll]_{\lambda^*}\\
%}
\label{eqn:ding105}
\end{equation}
is commutative.

From \eqref{eqn:ding81} it follows that
\[
\lambda^*(u) = \sum u^{q^k} \otimes \xi_k,
\]
which when raised to the $q^r$th power gives
\[
\lambda^*(u^r) = \sum u^{q^{k+r}} \otimes \xi_k^{q^r} ,
\]
and leads to the formula
\[
(\lambda^* \otimes 1)\big(\lambda^*(u)\big) =
	\big(\lambda^* \otimes 1\big)\Big(\sum_k u^{q^k} \otimes \xi_k\Big)
= \sum_r\sum_k u^{q^{k + r}} \otimes \xi_r^{q^k} \otimes \xi_k .
\]
Whereas, the other way around the diagram \eqref{eqn:ding105} yields
\[
\big(1 \otimes \nabla_*\big)(\lambda^*(u)) =
	\sum_j u^{q^j}\otimes \nabla_*(\xi_k) ,
\]
and equating these two expressions leads to the asserted formula
for the coproduct.
\end{proof}

As remarked at the end of the previous section
the operations $\P^{p^i}$ for $i > 0$
are indecomposables in $\PP^*$, so
$\PP^*$ is not generated by the operations $\P^{q^i}$ for $i \geq 0$;
we need all the $\P^{p^i}$ for $i > 0$.
This can be readily seen on hand from the dual Hopf
algebra, where, since $\F_q$ has characteristic $p$, the elements
$\xi_1^{p^i}$ for $i \geq 0$ are all primitive, \cite{MM}.
The following corollary also indicates that passing from the prime field
$\F_p$ to a general Galois field $\F_q$ is not just a simple substitution 
of $q$ for $p$.

\begin{corollary}
Let $p$ be a prime integer, 
$q=p^\nu$ a power of $p$s and
$\F_q$ the Galois field with $q$ elements. The indecomposable module 
$Q(\PP^*)$
of the Steenrod algebra of $\F_q$ has a basis consisting of the elements
$\P^{p^i}$ for $i \in \N_0$, and the primitive elements $P(\PP^*)$
has a basis
consisting of the elements $\P^{\rm\Delta_k}$ for $k \in \N$, where,
for $k \in \N$,
$\P^{\rm\Delta_k}$ is dual to $\xi_k$ with
respect to the monomial basis for $\PP_*$  .
\end{corollary}

\section{The Milnor basis and embedding one Steenrod algebra in another}
\label{sec:sec5}

%\sec{milbasis}
%% introduce the Milnor indexing scheme for the elements of the steenrod
%% algebra dual to the monomial basis
%% the "spreading out" construction then leads to the desired embedding
%% theorem of one steenrod algebra in the next

If $I = (i_1, i_2,\ldots , i_k, \dots)$ is an index sequence
we denote by $\P(I) \in \PP^*(\F_q)$ the element in the Steenrod algebra 
that
is dual to the corresponding monomial
$\xi^I$ in $\PP_*(\F_q)$ with respect to the monomial basis for
$\PP_*(\F_q)$. This is not the same as the monomial
$\P^I = \P^{i_1}\cdot\P^{i_2}\cdots\P^{i_k}\cdots$, these two elements 
do not
even have the same degrees. As $I$ ranges over all index sequences the
collection $\P(I)$ ranges over an $\F_q$--basis for $\PP^*(\F_q)$
called the \emph{ Milnor basis\/}.

To give some examples of elements written in the Milnor basis
introduce the index sequence $\Delta_k$ which has a $1$ in the $k$th
position and otherwise $0$s. Then $\P^k$ is $\P(k\cdot\Delta_1)$,
and, as noted at the end of \fullref{sec:sec4},
\emph{the Milnor primitive elements\/}
$\P^{\Delta_k}= \P(\Delta_k)$, for $k > 0$ ,
form a basis for the subspace of all primitive elements.
In terms of the reduced power operations these elements can also be 
defined by
the inductive formulae
\[
\P^{\Delta_k} = 
\begin{cases}
\P^1 & k = 1 \\
[\P^{q^{k-1}}, \P^{\Delta_k}] &  k > 0,
\end{cases}
\]
where $[\P',\P'']$ denotes the commutator
$\P'\cdot\P'' - \P''\cdot\P'$ of $\P'$ and $\P''$. In Milnor's paper one
can also find a formula for the product $\P(I)\cdot\P(J)$ of two elements 
in the
Milnor basis.
The basis transformation matrix from the admissible to the Milnor basis 
and its
inverse is quite complicated, so we will say nothing more about it.

To each index sequence $I$ we can make correspond
both an admissible sequence over $\F_p$ and one over $\F_q$
via the equations \eqref{eqn:ding67} from the previous section.
This correspondence gives us a map
$\theta \co \PP^* (\F_q) \lra \PP^* (\F_p) \otimes_{\F_p} \F_q$.

\begin{theorem}
Let $p$ be a prime integer, 
$q=p^\nu$ a power of $p$, and
$\F_q$ the Galois field with $q$ elements. The map
\[
\theta \co \PP^* (\F_q) \lra \PP^* (\F_p) \otimes_{\F_p} \F_q,
\]
embeds the Steenrod algebra
$\PP^*(\F_q)$ of $\F_q$ as a Hopf subalgebra in the Steenrod algebra of
$\F_p$ extended from $\F_p$ up to $\F_q$.
\label{thm:thm5.1}
\end{theorem}

\begin{proof}
It is much easier to verify that the dual map
\[
\theta_* \co \PP_* (\F_p) \otimes_{\F_p} \F_q \lra \PP_* (\F_q) ,
\]
which is defined by
the requirement that it be a map of algebras, and take the values
\[
\theta_* (\xi_k(p) \otimes 1) = 
\begin{cases}
\xi_m (q) & k = m\nu ~(\text{so } p^k - 1 = q^m - 1) \\
0 & \text{otherwise},
\end{cases}
\]
on algebra generators, is a map of Hopf
algebras. This is a routine computation.
\end{proof}

The Steenrod algebra over the prime field $\F_p$ has a 
well known
interpretation as the mod $p$ cohomology of the Eilenberg--MacLane 
spectrum.
By flat base change $\PP^*(\F_p)\otimes_{\F_p} \F_q$ may be regarded as the
$\F_q$--cohomology of the same. By including the Eilenberg--MacLane spectrum
$ K(\F_p)$ for the prime field into the Eilenberg--MacLane spectrum 
$ K(\F_q)$ we may view the elements of $\PP^*(\F_p)\otimes_{\F_p} 
\F_q$
as defining stable
cohomology operations in $\F_q$--cohomology. By \fullref{thm:thm5.1}
this also allows us to interpret elements of $\PP^*(\F_q)$ as stable
cohomology operations acting on the $\F_q$--cohomology of a topological 
space.
Which elements appear in this way is described in cohomological terms in
\cite{autos}.

\section{Closing comments}
\label{sec:sec6}
%\sec{close}

Algebraic topologists will of course immediately say {\sl ''but that isn't 
the
Steenrod algebra, it is only the algebra of reduced power operations; 
there is
no Bockstein operator unless $q=2$.''}
This is correct, the full Steenrod algebra, with the
Bockstein, has not yet played a significant role in invariant theory, so it
has not been treated here. But, if
one wishes to have a definition of the full Steenrod algebra in the same 
style
as the one presented here, all one needs to do for $q \neq 2$
is to replace the functor
$V \leadsto \F[V]$ with the functor
$V \leadsto H(V)$, where $H(V)$ is defined to be
$H(V) = \F[V] \otimes E[V]$, with $E[V]$ the exterior algebra on the dual
vector space $V^*$ of $V$. Since $V^*$ occurs {\it twice} as a subspace 
of
$H(V)$, once as $V^*\otimes \F \subset \F[V]\otimes \F$
and once as $\F \otimes V^* \subset E[V]$, we need a way
to distinguish these two copies. One way to do this is to write $z$ for a
linear form $z \in V^*$ when it is to be regarded as a polynomial 
function, and
$dz$ for the same linear form when it is to be regarded as an alternating
linear form. This amounts to identifying $H(V)$ with the algebra of 
polynomial
differential forms on $V$.

Next introduce the Bockstein operator
$\beta \co H(V) \lra H(V)$ by requiring it to be the unique derivation with 
the
property that
for an alternating linear form $dz$ one has
$\beta(dz) = z$, where $z$ is the corresponding polynomial linear form,
and for any polynomial linear form $z$ one has
$\beta(z) = 0$. The operators $\P^k$
for $k \in \N_0$ together with $\beta$ generate a subalgebra of the 
algebra of
\
endomorphisms of the functor $V \leadsto H(V)$, and this subalgebra
is the full Steenrod algebra of the Galois field $\F_q$.

Finally, 
%% at the summer school T.P. Bisson spoke about his work with A. Joyal on
there is a universal algebra approach to both the Dyer--Lashof algebra and 
the Steenrod
algebra in \cite{bissonjoyal}. The interested reader should consult this 
paper
which contains many informative facts.

\bibliographystyle{gtart}
\bibliography{link}
\end{document}